\documentclass{article}
\usepackage{amsfonts,amssymb}

\makeatletter
\renewcommand{\@biblabel}[1]{#1.} 
\makeatother

\newcommand{\sgn}{\mathop{\mathrm{sgn}}}
\newcommand{\dsp}{\displaystyle}
\newcommand{\qed}{\hfill $\square$}

\newtheorem{theorem}{Theorem}[section]
\newtheorem{lemma}[theorem]{Lemma}

\newtheorem{corollary}[theorem]{Corollary}
\newtheorem{remark}[theorem]{Remark}

\newcommand{\R}{{\mathbb R}}

\title{\bf   
Spherical solutions to the Klein-Gordon equation in the expanding universe
}

\author{{\bf Karen Yagdjian} 
 }

\begin{document}

\date{}
\maketitle
\thispagestyle{empty}
\vspace{-0.3cm}

\begin{center}
{\it School of Mathematical and Statistical Sciences,\\
University of Texas RGV,
1201 W.~University Drive, \\
Edinburg, TX 78539,
USA }
\end{center}
\medskip

\renewcommand{\theequation}{\thesection.\arabic{equation}}
\setcounter{equation}{0}
\pagenumbering{arabic}
\setcounter{page}{1}
\thispagestyle{empty}

\hspace{2cm}\begin{abstract}
\begin{small}

\medskip

We produce an explicit formula for the wave function of the spherically symmetric fields emitted to the FLRW universe with the scale factor generated by the de~Sitter universe. As an application of these explicitly written solutions of the Klein-Gordon equation, we   test   the decay in time  of the field  generated  by a pionic atom.
\medskip

\medskip

\end{small}
\end{abstract}

\setcounter{equation}{0}

\section{Introduction. Main Results}

In this study, we produce an explicit wave function of the particles (fields) emitted to the  Friedmann-Lema\^itre-Robertson-Walker (FLRW)  universe with the scale factor generated by the de~Sitter space. We consider models of  scalar fields with spherical symmetry at the initial time.  They can be generated, for instance, by a pionic atom,  exotic atom,  or  gravitational field (see, e.g., \cite{Costa-Alho-Natario,Ericson,Nikiforov-Uvarov,JMP_2024}). 
  In the expanding universe, the particles obey the equation for the free particles. The purpose of this article is to derive the explicit form of exact solutions of the spherical waves in curved spacetime. The models illustrate a transition from the initial spacetime to the non-static FLRW universe with a de~Sitter scale factor.

Thus, for the scalar field,  we are interested in the Cauchy problem 
\[ 
\cases{ \dsp \Box_g \Phi - \frac{m_n^2 c^2}{\hbar^2}  \Phi=  0\,, \cr
 \dsp \Phi (x,0)=  F(r,0)  Y_{\ell m}(\theta, \phi)\, , \quad \Phi _t(x,0)=F_t(r,0)  Y_{\ell m}(\theta, \phi)\,,}
\]
where $m_n$ is a mass, $F (r,t) $ is a time-dependent radial part of the initial wave, $F_t(r,t):={\partial_t F }(r,t)$, while 
$ Y_{\ell m}(\theta, \phi)$ are spherical harmonics. The line element of the non-stationary spacetime  with the metric tensor $ g_{ik}$ and the scale factor $a=a(t) $ in the spherical coordinates is
\begin{equation}
 \label{ds_NM}
\dsp ds^2= -  c^2dt^2
+  a^2(t) dr^2 + a^2(t) r^2(d\theta ^2 + \sin^2 \theta \, d\phi ^2) .
\end{equation}
For the de~Sitter  universe in FLRW form, $a(t)= e^{Ht}$, where $H$ is the Hubble parameter. The scalar curvature of this space is $12H^2$.   Consequently, the covariant d'Alembert's operator in the spacetime with the metric  of (\ref{ds_NM}) is 
\[
\square_g  \Phi = \frac{1}{\sqrt{|g|}}\frac{\partial }{\partial x^\alpha}\left( \sqrt{|g|} g^{\alpha \beta} \frac{\partial  \Phi }{\partial x^\beta} \right),
\]
where $ \alpha$ and $\beta$ are running  from $0$ to $3$. We are interested in the  explicit form of the spherical field in the de~Sitter FLRW spacetime, which is  a solution  of the  problem 
\begin{equation}
\label{1.10n}
\cases{ \dsp  \Phi_{tt} +   3  H \Phi_t - c^2e^{-2H t} \Delta \Phi +  \frac{m_n^2 c^4}{\hbar^2}\Phi=  0\,, \cr
\dsp  \Phi (x,0)= F_0 (r)Y_{\ell m}(\theta, \phi)\, , \quad \Phi _t(x,0)=  F_1(r) Y_{\ell m}(\theta, \phi) \,.}
\end{equation}
Here, $Y_{ lm} (\theta, {\phi})$ is a spherical harmonic, $\ell =0,1,\ldots$, $ m=0,\pm 1,\pm 2 \ldots$, $|m|\leq \ell$. In particular, the initial data $\Phi (x,0) $ and $\Phi _t(x,0) $ can be given by a wave function 
of  the pionic  atom in Minkowski space, that is, by the functions $
  \frac{R_{n\ell}(r)}{r} Y_{\ell m}(\theta, \phi)
$ and $
  \left( -i\frac{E_n  }{\hbar}\right)  \frac{R_{n\ell}(r)}{r} Y_{\ell m}(\theta, \phi)
$, respectively  (see, e.g., \cite{Ericson}). We will use a system of units in which the speed of light $c=1$.  

Denote $\dsp M :=\left(\frac{9}{4}H^2 -\frac{m_n^2 }{\hbar^2} \right)^{1/2}$. 
Let $F(a,b;c;z) $ be a hypergeometric function (see, e.g., \cite{B-E}), and $ J_{\nu}(z)$ be a Bessel function. We define  
 kernels   $K_0(z,t;M) $   and  $K_1(z,t;M) $   as~\footnote{There is a typo in  $K_0$ on page 11~\cite{JMP_2024}; the correct version is in Lemma~3.1~\cite{JMP_2024} }
\begin{eqnarray}
\label{K0M}
K_0(r,t;M)
&  := &
-4^{-\frac{M}{H}} e^{-H t} e^{M t} \left(\left(e^{-H t}+1\right)^2-H^2 r^2\right)^{\frac{M}{H}-\frac{5}{2}}\\
&  &
\times \Bigg\{ e^{H t} \left(\left(e^{-H t}+1\right)^2-H^2 r^2\right) \left(-H^2 M r^2+M e^{-2 H t}+H e^{-H t}+H-M\right) \nonumber\\
&  &
\times  F \left(\frac{1}{2}-\frac{M}{H},\frac{1}{2}-\frac{M}{H};1;\frac{\left( 1-e^{-H t}\right)^2-H^2 r^2}{\left(1+e^{-H t}\right)^2-H^2 r^2}\right)\nonumber\\
&  &
+\frac{1}{H}(H-2 M)^2 \left(-H^2 r^2+e^{-2 H t}-1\right)\nonumber\\
&  &
\times F\left(\frac{3}{2}-\frac{M}{H},\frac{3}{2}-\frac{M}{H};2;\frac{\left(1-e^{-H t}\right)^2-H^2 r^2}{\left(1+e^{-H t}\right)^2-H^2 r^2}\right)\Bigg\},\nonumber\\
\label{K1M}
K_1(r,t,M)
& :=  &
4^{-\frac{M}{H}} e^{M t} \left(\left(e^{-H t}+1\right)^2-H^2 r^2\right)^{\frac{M}{H}-\frac{1}{2}}\\
&  &
\times   F \left(\frac{1}{2}-\frac{M}{H},\frac{1}{2}-\frac{M}{H};1;\frac{\left(1- e^{-H t}\right)^2-H^2 r^2}{\left(1+e^{-H t}\right)^2-H^2 r^2}\right)\,.  \nonumber
\end{eqnarray}
We  often omit the parameter $M$ in  $K_0(z,t;M) $   and  $K_1(z,t;M) $ when it  does not lead to confusion. 
Denote $\phi (t):= (1-e^{-Ht})/H$ a conformal time. Note: the scalar curvature $R$ is encoded in these kernels via the Hubble constant $H=\sqrt{R/12} $. The main results of this study, which are obtained by the Integral Transform Approach (ITA) \cite{Rend_2010,MN2015}, are the following two theorems.

\begin{theorem}
\label{T1.1} Assume that $ F_0(r),  F_1(r) \in    C^1  ({\mathbb R}\setminus \{0\})$, $\ell$ is non-negative  integer number, and the condition 
\[
F_i(r)=O(r^{\mu -\frac{1}{2}}) \,\, \mbox{\rm as } \,\, r \searrow 0 ,\quad i=0,1\,,
\] 
 is fulfilled  with $\mu  >\ell-\frac{3}{2}$. 
Further,   assume that   $ r^ \ell   F_0(r)$ and $r^ \ell F_1(r) $  are     even functions. 
If   $Y_{ lm} (\theta, {\phi})$ is a spherical harmonic, $\ell=0,1,\ldots$, $ m=0,\pm 1\pm 2,\ldots$, and $|m|\leq \ell$,  then the spherical field in the de~Sitter FLRW spacetime  
is given as follows
\begin{eqnarray*} 
&  &
\Phi (r,\theta, \phi,t) \\ 
& = &
e^{-Ht} Y_{\ell m}(\theta, \phi)\Bigg\{ \frac{1}{2r} [ (r-\phi (t) )F_0(r-\phi (t) )   +  (r+\phi (t))F_0(r+\phi (t)) ]\\
&  &
-\frac{1}{4 } \ell  (\ell +1) \frac{\phi (t)}{r^2}\int_{r-\phi (t)}^{r+\phi (t)}  F_0(s)  
F \left(1-\ell ,\ell +2;2;\frac{\phi (t)^2-(r-s )^2}{4 r s }\right)ds\Bigg\} \\
&  &
+ \,  e^{-H\frac{3}{2}t}Y_{\ell m}(\theta, \phi)\int_{ 0}^{\phi (t)}  \Bigg\{ \frac{1}{2r} [ (r-\tau )F_0(r-\tau )   +  (r+\tau)F_0(r+\tau) ]\\
&  &
-\frac{1}{4 } \ell  (\ell +1) \frac{\tau}{r^2}\int_{r-\tau}^{r+\tau}  F_0(s)  
F \left(1-\ell ,\ell +2;2;\frac{t^2-(r-s )^2}{4 r s }\right)ds\Bigg\} \left(2  K_0 \left(\tau ,t \right)+ 3K_1\left(\tau ,t \right)\right)\,  d\tau  \nonumber \\
& &
+\, 2e^{-\frac{3}{2}Ht}Y_{\ell m}(\theta, \phi)\int_{0}^{\phi (t)}   \Bigg\{ \frac{1}{2r} [ (r-\tau )F_1(r-\tau )   +  (r+\tau)F_1(r+\tau) ]\\
&  &
-\frac{1}{4 } \ell  (\ell +1) \frac{\tau}{r^2}\int_{r-\tau}^{r+\tau}  F_1(s)  
F \left(1-\ell ,\ell +2;2;\frac{\tau^2-(r-s )^2}{4 r s }\right)ds\Bigg\}
  K_1\left(\tau,t \right)\, d\tau\,.
\end{eqnarray*}

For the huygensian scalar field (particle) (\ref{1.10n}), that is, $ m_n^2   =2H^2\hbar^2 $,  
the spherical field in the de~Sitter FLRW spacetime is 
\begin{eqnarray*} 
&  &
\Phi  (r,\theta, \phi,t) \\
& = &
e^{-Ht}  Y_{\ell m}(\theta, \phi)\Bigg\{ \frac{1}{2r} [ (r-\phi (t) )F_0(r-\phi (t) )   +  (r+\phi (t))F_0(r+\phi (t)) ]\\
&  &
-\frac{1}{4 } \ell  (\ell +1) \frac{\phi (t)}{r^2}\int_{r-\phi (t)}^{r+\phi (t)}  F_0(s)  
F \left(1-\ell ,\ell +2;2;\frac{\phi (t)^2-(r-s )^2}{4 r s }\right)ds\Bigg\}\\
&  &
+ \,   e^{-Ht}   
 Y_{\ell m}(\theta, \phi)\int_0^{\phi (t)} \Bigg\{ \frac{1}{2r} \Big( (r-\tau )[F_0(r-\tau ) +F_1(r-\tau )]  +  (r+\tau)[F_0(r+\tau)+F_1(r+\tau)] \Big)\\
&  &
-\frac{1}{4 } \ell  (\ell +1) \frac{\tau}{r^2}\int_{r-\tau}^{r+\tau}  (F_0(s)+F_1(s))   
F \left(1-\ell ,\ell +2;2;\frac{\tau^2-(r-s )^2}{4 r s }\right)ds\Bigg\}d \tau \,.
\end{eqnarray*}
\end{theorem}
 In the case of $ m_n^2   =2H^2\hbar^2 $, the Klein-Gordon equation
becomes the conformally invariant wave equation.  
It should be noted here that for negative integers $a$, the hypergeometric function $F(a,b;c;z) $ is a polynomial of degree $-a$ (see, e.g., \cite{B-E}). 

In the next theorem, by imposing conditions on the behavior of the functions $F_0= F_0(r)$ and $F_1= F_1(r) $ for large $r$ one can appeal to the Hankel transform and write the representation of the field in the form that is more amenable for the $L^2({\mathbb R}^3)$-estimates. In particular, this applies to the wave function of the pionic atom that  arrives in the de~Sitter FLRW spacetime, because that wave function has an exponential decay at  infinity of spatial coordinates (see, e.g., \cite{Ericson}).       
\begin{theorem}
\label{T3.1}
Assume that $ F_0(r),  F_1(r) \in    C^1  (0,\infty)$, and 
\[
\int_0^\infty  r ^{3/2}|F_i(r)|\,dr<\infty,\quad    F_i(r)=O(r^{-k}), \quad r \to \infty, \quad k> 2,\quad i=0,1.
\] 
If   $Y_{ lm} (\theta, {\phi})$ is a spherical harmonic, $\ell=0,1,\ldots$, $|m|\leq \ell$,  then the spherical field in the de~Sitter FLRW spacetime with the initial state
\[
\Phi (x,0)= F_0(r) Y_{\ell m}(\theta, \phi)\, , \quad
\Phi _t(x,0)=F_1(r)Y_{\ell m}(\theta, \phi)
\]
is given by
\begin{eqnarray*} 
&  &
\Phi (r,\theta, \phi,t) \\  
& = &
e^{-Ht} Y_{\ell m}(\theta, \phi)\frac{1}{\sqrt{r}}\int_{0}^\infty \left(\int_{0}^\infty      F_0(\rho) J_{\ell+\frac{1}{2}}(\rho \lambda )\rho^{3/2}\, d\rho\right)\cos \left( \lambda \phi (t)  \right) J_{\ell+\frac{1}{2}}(r \lambda )\lambda\, d\lambda \\
& &
+ \,  e^{-\frac{3}{2}Ht}Y_{\ell m}(\theta, \phi)\frac{1}{\sqrt{r}}\int_{ 0}^{1}\int_{0}^\infty \left(\int_{0}^\infty      F_0(\rho) J_{\ell+\frac{1}{2}}(\rho \lambda )\rho^{3/2}\, d\rho\right)\cos \left( \lambda \phi (t)  s \right) \\
&  &
\times  J_{\ell+\frac{1}{2}}(r \lambda )\lambda\, d\lambda\left(2  K_0 \left(\phi (t) s,t \right)+ 3K_1\left(\phi (t) s,t \right)\right)\phi (t) \,  ds  \nonumber \\
& &
+\, 2e^{-\frac{3}{2}Ht}Y_{\ell m}(\theta, \phi)\frac{1}{\sqrt{r}}\int_{0}^1 \int_{0}^\infty \left(\int_{0}^\infty      F_1(\rho) J_{\ell+\frac{1}{2}}(\rho \lambda )\rho^{3/2}\, d\rho\right)\cos \left( \lambda \phi (t)  s \right)  \\
&  &
\times 
 J_{\ell+\frac{1}{2}}(r \lambda )\lambda\, d\lambda \, K_1\left(\phi (t) s,t \right)\phi (t) \, ds \quad \mbox{ for all} \quad x \in {\mathbb R}^3, \,\, t>0 . 
\end{eqnarray*}
For the huygensian scalar field (particle) (\ref{1.10n}), that is, $ m_n^2   =2H^2\hbar^2 $,  
the spherical field in the de~Sitter FLRW spacetime is 
\begin{eqnarray*}
&  &
\Phi  (r,\theta, \phi,t)\\
  & = &
e^{-Ht}  Y_{\ell m}(\theta, \phi)\frac{1}{\sqrt{r}}\Bigg\{ \int_{0}^\infty \left(\int_{0}^\infty      F_0(\rho) J_{\ell+\frac{1}{2}}(\rho \lambda )\rho^{3/2}\, d\rho\right)\cos \left( \lambda \phi (t)  \right) J_{\ell+\frac{1}{2}}(r \lambda )\lambda\, d\lambda \\
&  &
+ 
\int_0^\infty 
  \left( \int_0^\infty (   F_0(\rho)+F_1(\rho))  J_{  \ell+\frac{1}{2}}(\lambda  \rho  ) \,\rho ^{3/2} d\rho  \right) \sin \left(\lambda \phi (t)  \right)   J_{ \ell+\frac{1}{2}}(r \lambda )\,d\lambda \Bigg\} \,. 
\end{eqnarray*}
\end{theorem}
\medskip

In Corollary~\ref{C3.2},   we prove the exponential decay of the solution's tail in the case where the radial part of the initial function is generated by the wave function of the pionic atom. Moreover, it is also shown that the obtained decay is optimal.  
\medskip

Below, we provide some examples of such an initial state $(\Phi (x,0), \Phi _t(x,0))$ of the field. In particular, it can be a wave function of a pionic atom in the Minkowski spacetime (see, e.g., \cite{Ericson,Greiner}), or a scalar field in the matter- or radiation-dominated or de~Sitter universe with a different Hubble constant $H_0\not= H$ and different mass. In  exotic atoms,   one of the electrons has been replaced by a negatively charged particle such as a muon, K-meson, pion, $\Sigma$-hyperon, etc. For further details on the physical interest in such systems, see, for example, \cite{Ericson,Greiner,G-S,Nancy}. A pionic atom is a hydrogen-like system  in which  an electron replaced by a spinless relativistic particle pion. A pion is an example of a particle with electromagnetic interactions that obeys the Klein-Gordon equation. In fact, by pionic atom, one can test the validity of the Klein-Gordon equation (see, e.g., \cite[Sec.18.4]{Greiner}). In (\ref{1.10n}), the pion has  left the atom, the Coulomb field is no longer in action, and spacetime has  become inflationary. Thus, the free pion propagates in the expanding universe  and  carries information about its previous state. 
\medskip

Although  spherically symmetric solutions for the Klein-Gordon and Dirac equations in the Minkowski spacetime have been known for a long time (see, e.g., \cite[\S 24]{Berestetskii-Lifshits}, \cite{deRham}, \cite{Ovsiyuk}, and references therein), explicit exact formulas in the curved spacetimes are known only for a few   special cases (see, e.g., \cite{Barut-D,deOliveira-Schmidt,Villalba2002,AnPH2020}). We believe that the obtained explicit formulas of Theorem~\ref{T3.1} can be useful, in particular, in the study of caustics for spherical waves in curved spacetime (see, e.g., \cite{deRham}). In the forthcoming paper, the representation formulas of Theorems~\ref{T1.1}, \ref{T3.1} will be used for the derivation of the exact spherical solutions to the spin-$\frac{1}{2}$ fields in the expanding universe.
\medskip

Another possible application of our formulas can be to the problem of the  existence of quasinormal modes (QNMs) on a non-static FLRW universe with a de~Sitter scale factor. 
The existence of  QNMs  on de Sitter   spacetime has been much studied in the literature; a detailed  analysis of the literature on the spherical waves and QNMs     has been   given in \cite{Hintz-Xie}.  
In \cite{Hintz-Xie}, it has been  demonstrated that QNMs  exist for massive
scalar fields on static de Sitter space in all spacetime dimensions greater than one and for all  complex-valued scalar field masses. The Taylor expansion    with the recursion relation for the coefficients has been presented for calculating of QNMs and their corresponding mode solutions. The proof of the convergence of the series has been done for the analytic data.    
\medskip

The   study of the spherically symmetric solutions of the wave equation in the de Sitter space in the Bondi coordinates  was carried out  in \cite{Costa-Alho-Natario}. In those coordinates, the problem was reduced to an integro-differential equation with the initial data on a light cone.  The main results state  existence and uniqueness, as well as decay of the solution in the Bondi coordinates and its limit at infinity.  
\medskip

The rest of this paper is organized as follows. In Section~\ref{SWM}, we describe the explicit formulas for spherical solutions to the wave equation in the Minkowski space. We list three different approaches in subsections~\ref{SS2.1}, \ref{SS2.2}, and \ref{SS2.3}, more exactly, approaches by the general solution, by the Riemann function, and by the Hankel transform, respectively. In Section~\ref{S3} we apply the ITA to write the explicit form of the Klein-Gordon fields in the Minkowski and de~Sitter spacetimes. This completes the proof of Theorems~\ref{T1.1},\ref{T3.1}.

\section{Spherical solutions of the wave equation in  \\ Minkowski spacetime}
\label{SWM}
\setcounter{equation}{0}

To make the presentation of the proof of Theorem~\ref{T3.1} self-contained, in this section,  we write explicit solutions to the wave equation with the initial data, which are the products of the spherical harmonic and the radial function. More information on this type of result can be found, for instance, in the review \cite{Stewart}, or in \cite{Airapetyan} and \cite{Witt}.    In our study,  the radial function can be a radial part of the wave function of the free spinless particle, or  pionic atom or can be generated by another model.  

Consider the problem for the wave equation
\begin{eqnarray} 
\label{WE}
V_{tt}-  \Delta V =0, \quad V(x,0)=V_0(x) , \quad V_t(x,0)=V_1(x)\,.
\end{eqnarray}
Here $ \Delta $ is the Laplace operator in spherical coordinates 
\[
\Delta=  \partial_r^2 +\frac{2     }{r}\partial_r+ \frac{  1 }{r^2}\partial_{\theta }^2
+\frac{\cot (\theta )    }{r^2}\partial_\theta+\frac{\csc ^2(\theta )   }{r^2} \partial_\phi^2 =  \partial_r^2 +\frac{2     }{r}\partial_r
+\frac{1}{r^2}\Delta_{\theta, \phi}\,.
\] 
The initial functions are assumed to be of the form 
\begin{equation}
\label{IVP_MSP}
V_0(x)=F_0(r)Y_{\ell m}(\theta, \phi),\quad V_1(x)=F_1(r)Y_{\ell m}(\theta, \phi),
\end{equation}
where $ Y_{\ell m}(\theta, \phi)$, $| m|\leq \ell$, $\ell=0,1,\ldots$,  are spherical harmonics,
\begin{eqnarray*}  
&  &
\Delta_{\theta, \phi} Y+\ell(\ell+1) Y=0,\quad \int\left|Y_{\ell m}(\theta, \phi)\right|^2 d \Omega=1.
\end{eqnarray*}
The functions $F_0$ and $F_1$ are smooth in $\R_+$ and  satisfy   certain  conditions as $r\to 0$ and $r\to \infty$. 
Consequently, we look for a solution of the form
\[ 
V(x,t)=F(r,t)Y_{\ell m}(\theta, \phi)\,.
\]
Then
\[
\cases{\dsp Y_{\ell m}(\theta, \phi)\frac{\partial^2}{\partial t^2}F(r,t)
-     Y_{\ell m}(\theta, \phi)\frac{1}{r}\frac{\partial^2}{\partial r^2}(rF(r,t))-F(r,t)\frac{1}{r^2}\Delta_{\theta, \phi}Y_{\ell m}(\theta, \phi)   =0, \cr 
F(r,0)= F_0(r), \quad F_t(r,0)= F_1(r), \quad r \geq 0\,.}  
\]
Hence, 
\begin{eqnarray} 
\label{EqF} 
&  &
 F_{tt}(r,t)-F_{rr}(r,t)
-    \frac{2}{r} F_r(r,t)+\frac{\ell(\ell+1)}{r^2}F(r,t)   =0\,. 
\end{eqnarray}
For the pionic atom the wave function is 
\[
\psi(x,t)=e^{itE}F(r,t)Y_{\ell m}(\theta, \phi)\,,
\]
where for the function $F=F (r,t)$ the radial equation is 
 \begin{eqnarray*}  
&  &
 F_{tt}(r,t)-F_{rr}(r,t)
-    \frac{2}{r} F_r(r,t))+\left( U(r) +\frac{\ell(\ell+1)}{r^2} \right) F(r,t)   =0,\quad r > 0,\quad t \in\R\,, 
\end{eqnarray*}
while $ U$ is a potential. For the function $R(r,t)=rF(r,t)$, the equation is 
\begin{eqnarray*}  
&  &
 R_{tt}(r,t)-R_{rr}(r,t)
+ \left( U(r) +\frac{\ell(\ell+1)}{r^2} \right) R(r,t)   =0,\quad r > 0\,. 
\end{eqnarray*}
The normalization conditions   are    
\begin{equation}
\label{integrable}
\int_0^\infty R^2(r,t) \,dr <\infty,\qquad  
\int_0^\infty F^2(r,t)r^2 \,dr <\infty. 
\end{equation}
 
\subsection{The general solution approach}
\label{SS2.1}

\bigskip

Consider  for the wave equation (\ref{WE})  the problem with  the initial conditions
\begin{eqnarray*} 
V(x,0)=F_0(r) Y_{\ell m}(\theta, \phi) , \quad V_t(x,0)=F_1(r) Y_{\ell m}(\theta, \phi)\,.
\end{eqnarray*}
We look for the solution $F=F(r,t) $ of equation (\ref{EqF})  satisfying (\ref{integrable}). For example, if  $F_0(r) $ is generated by the wave function of the pionic atom \cite[pp.188-189]{Ericson}, \cite[Sec.18.4]{Greiner}, and \cite[ \textsection 26, Sec 2]{Nikiforov-Uvarov} in the Minkowski space, then
\begin{eqnarray}
\label{aroundzeroEric_Exp}
\hspace{-0.7cm}&  &
F_0(r)= C_{n,\ell,\mu} r^{\mu -\frac{1}{2} } e^{-r/2 }L_{n-\ell-1}^{2\mu }(r), \quad L_{n-\ell-1}^{2\mu }(r)\quad \mbox{\rm is Laguerre polynomial},\,\, r>0 ,
\end{eqnarray}
where with the fine structure constant $\alpha = e^2/(\hbar c) \approx  1/137$ and  the charge $e$ of the electron in Gaussian units, 
\begin{eqnarray*}
\hspace{-0.7cm} &  & 
\mu= \sqrt{\left(\ell+\frac{1}{2}\right)^2- Z^2 \alpha^2} \geq 0\,. 
\end{eqnarray*} 
If $\ell=0$, while  the atomic number $Z=1$ and $Z=2$,    then $\mu-\frac{1}{2}=-0.0000532822$ and $\mu-\frac{1}{2}=-0.000213163$, respectively.

For the case of  $F_1(r)=0 $, assume that $F_0(r) $ is  smooth on  $\R \setminus\{0\}$.   
Define the operator
\begin{eqnarray*} 
L_\ell (r,\partial_r,\partial_t) & := &
  \partial_t ^2- \partial_r ^2
-    \frac{2}{r} \partial_r+\frac{\ell(\ell+1)}{r^2} \,. 
\end{eqnarray*}
It is easily seen that
\begin{eqnarray*} 
 \left(\partial_r-\frac{\ell }{r}\right)L_{\ell }
& = &
L_{\ell +1} \left(\partial_r-\frac{\ell }{r}\right)\,. 
\end{eqnarray*}
Hence, 
if $w_{(\ell)} =w_{(\ell)}(r,t) $ is a solution to $L_\ell (r,\partial_r,\partial_t)w_{(\ell)}(r,t)=0 $, then the function 
\[
w_{(\ell+1)}(r,t)=\left( \partial_r-\frac{\ell}{r}\right)w_{(\ell)}(r,t)
\]   
solves the equation 
\[
L_{\ell+1} (r,\partial_r,\partial_t)w_{(\ell+1)}(r,t)=0 . 
\]
Therefore,  
the general solution to equation (\ref{EqF}) is given by the formula
\begin{eqnarray*} 
&  &
F ^{(\ell)}(r,t)= r^\ell\left( \frac{1}{r} \frac{\partial}{\partial r}\right)^\ell\left( \frac{\Phi (r+t)+\Psi(r-t)}{r}\right),  \quad r \in {\mathbb R},  \quad t \geq 0\,, \quad \ell=0,1,2,\ldots\,,
\end{eqnarray*}
where $ \Phi$ and $\Psi $ are the arbitrary functions. For every given $\ell$,  one can derive from the general solution an explicit formula for the solution to the Cauchy problem
(\ref{WE})\&(\ref{IVP_MSP}), provided that (\ref{integrable}) holds. The calculations lead to the following result. 
The solution $F^{(\ell)}(r,t) $, $\ell =0,1,2,\ldots$,  of  problem (\ref{EqF}) with the initial data
\[
F^{(\ell)}(r,0)=F_0^{(\ell)}(r),\qquad F_t^{(\ell)}(r,0)=0,\qquad r \in {\mathbb R}_+\,,
\]
is given recursively by 
\begin{equation}
\label{Fellplus1}
F^{(\ell+1)}(r,t)=\left( \partial_r-\frac{\ell}{r}\right)\widetilde{F}^{(\ell )}(r,t),\qquad \ell =0, 1,2,\ldots\,,
\end{equation}
where $\widetilde{F}^{(\ell )}(r,t)$ is given by  
\begin{eqnarray}
\label{Felltilde}  
\widetilde{F}^{(\ell )}(r,t)
& = &
\frac{1}{2 r }\left\{(r-t) \widetilde{F}_0^{(\ell )}(r-t) +(r+t) \widetilde{F}_0^{(\ell )}(r+t)\right\} \\
&  &
+
\frac{1}{2 r^{\ell+1}}t \cases{ \dsp \sum_{   k \in {\mathbb Z} \,\,is\,\, odd ,\,\, |k|\leq \ell-1 \atop{\alpha+\beta+k=\ell-1 \atop{\alpha \geq 0,\beta \geq 0\,\,are\, \,  even}} \,\,}a_{ k,\alpha,\beta}^{(\ell )}r^{\alpha }t^{\beta} \int _{r-t}^{r+t}  \tau ^{k}\widetilde{F}_0^{(\ell )}(\tau ) \,d\tau  ,\quad \ell \,\, \mbox{is even},\cr
{}\cr
\dsp \sum_{ \alpha+\beta+k=\ell-1\atop{\atop{k \in {\mathbb Z} \,\,is\,\, even ,\,\, |k|\leq \ell-1 \atop{\alpha \geq 0,\beta \geq 0 \,\,are\, even}}  } \,\,}b_{ k,\alpha,\beta}^{(\ell )}r^{\alpha }t^{\beta} \int _{r-t}^{r+t} \tau ^{k}\widetilde{F}_0^{(\ell )}(\tau )\,d\tau,\quad \ell \,\, \mbox{is odd}} \nonumber
\end{eqnarray}
and
\begin{equation}
\label{FtileF}
\widetilde{F}_0^{(\ell )}(r)
=r^\ell \int _0^rs^{-\ell}F_0^{(\ell+1 )}(s)\, ds\,.
\end{equation}
Here $a_{ k,\alpha,\beta}^{(\ell )} $ and $b_{ k,\alpha,\beta}^{(\ell )} $ are some suitable coefficients.

 For   the first three values of $\ell$, if $F_0 $ is  given by (\ref{aroundzeroEric_Exp}), then the function $F (r,t)$ can be a distributional solution from the Sobolev space. For the first five values of $\ell$, the results  are as follows. 
For  $\ell=0$ and the even function  $ F_0(r)$,  the  solution is 
\begin{eqnarray*}
F (r,t)
& =  &   
 \frac{1}{2r} \left\{(r-t) F_0(r-t)+  (r+t) F_0(r+t)  \right\} ,  \quad t \geq 0\,.
\end{eqnarray*}
  For $\ell=1$ and  odd function $ F_0(r)$,  the  solution is  
\begin{eqnarray*}
F(r,t)
&  =& 
\frac{1}{2 r }\left\{   (r-t) F_0(r-t)+  (r+t) F_0(r+t)\right\}  -\frac{1}{2 r^2}t  \int_{r-t} ^{r+t} F_0(\tau ) \, d\tau  \,.
\end{eqnarray*}
For $\ell=2$ and even function $ F_0(r)$,  the  solution is 
\begin{eqnarray*}
F(r,t)
& = &
\frac{1}{2 r }\left\{  (r-t) F_0(r-t)+ (r+t) F_0(r+t)\right\} \nonumber \\
&  &
+\frac{1}{2 r^3}t\Bigg\{ -\frac{3}{2}(r^2-t^2) \int_{r-t}^{r+t} \frac{F_0(\tau )}{ \tau } \, d\tau 
 -  \frac{3}{2}  \int_{r-t}^{r+t} \tau  F_0(\tau ) \, d\tau \Bigg\} \,.
\end{eqnarray*}
For $\ell=3$ and odd function $ F_0(r)$,  the  solution is 
\begin{eqnarray*}
F(r,t)
& = &
 \frac{1}{2 r }\left\{ (r-t) F_0(r-t)+  (r+t) F_0(r+t)\right\} \nonumber\\
 &  &
\frac{1}{2 r^4} t\left\{ -\frac{15}{8} (r^2-t^2)^2  \int _{r-t}^{r+t}\frac{F_0(\tau )}{  \tau ^2}d\tau 
  \nonumber \right. \\
&  &
\left. -\frac{1}{4}  \left(9 r^2 -15 t^2\right) \int _{r-t}^{r+t} F_0(\tau )d\tau 
 -\frac{15}{8}  \int _{r-t}^{r+t} \tau ^2 F_0(\tau )d\tau \right\}\,.  
\end{eqnarray*}
For $\ell=4$ and even function $ F_0(r)$,  the  solution is  
\begin{eqnarray*}
F(r,t) 
& = &
\frac{1}{2 r }\left\{  (r-t) F_0(r-t) +   (r+t) F_0(r+t)\right\} \nonumber\\
&  &
-\frac{1}{2 r^5}t\left\{ \frac{35}{16 } \left(r^6-3 r^4 t^2+3 r^2 t^4-t^6\right)\int _{r-t}^{r+t}\frac{F_0(\tau )}{ \tau ^3}d\tau \nonumber \right. \\
&  &
+ \frac{5}{16} \left(9 r^4-30 r^2 t^2+21 t^4\right) \int _{r-t}^{r+t}\frac{F_0(\tau )}{\tau }d\tau \nonumber \\
&  &   
\left. +\frac{5}{16}\left(9 r^2-21 t^2\right)  \int _{r-t}^{r+t} \tau  F_0(\tau )d\tau +\frac{35 }{16} \int _{r-t}^{r+t}\tau ^3 F_0(\tau )d\tau   \right\}\,. \nonumber 
\end{eqnarray*}

The obtained formulas allow us to estimate the 
``radial tail''  of the solution.
\begin{lemma} 
\label{LTail}
Assume that $F_1(r)\equiv 0 $.  
For a pionic atom with (\ref{aroundzeroEric_Exp}), the tail undergoes an exponential decay in time 
\[ 
\left|F(r,t)- \frac{1}{2 r }\left\{   (r-t) F_0(r-t)+  (r+t) F_0(r+t)\right\} \right|
  \leq   
C(r)e^{-t/2}t^{n +\mu - \frac{3}{2} }   , \quad \ell=1,2,\ldots  \,. 
\] 
\end{lemma}
\medskip
 
 \noindent
 {\bf Proof.} 
Indeed,  for even $\ell$ the tail contains the integrals 
\begin{eqnarray*} 
 \int _{r-t}^{r+t} F_0(\tau )  \tau ^k d\tau 
& =  & 
\int _{t-r}^{t+r} F_0(\tau )  \tau ^k d\tau   , \quad k=-\ell+1,-\ell+3,\ldots,  \ell-3,\ell-1  \,. 
\end{eqnarray*}
For a pionic atom with (\ref{aroundzeroEric_Exp}) for every given $r$ and a sufficiently large time $t>r+1$,  
\begin{eqnarray*} 
\left| \int _{r-t}^{r+t} F_0(\tau )  \tau ^k d\tau \right| 
& \leq   & 
\int _{t-r}^{t+r} \left| F_0(\tau ) \tau ^k \right| d\tau  \\
& \leq   & 
C_{n,\ell,\mu}\int _{t-r}^{t+r} \left| \tau^{\mu -\frac{1}{2} } e^{-\tau/2 }L_{n-\ell-1}^{2\mu }(\tau) \tau ^k \right| d\tau  \\
& \leq   & 
C_{n,\ell,\mu}\int _{t-r}^{t+r}  e^{-\tau/2 }  \tau^{a}    d\tau \,, 
\end{eqnarray*}
where $ a:=\mu -\frac{1}{2}+ n-\ell-1+k$. The last integral can be evaluated in terms of the   incomplete gamma function (\cite[Ch.9]{B-E}): 
\begin{eqnarray*}
\left| \int _{r-t}^{r+t} F_0(\tau )  \tau ^k d\tau \right| 
& \leq   & 
C_{n,\ell,\mu}2^{a+1} \left|\Gamma \left(a+1,\frac{t-r}{2}\right)-\Gamma \left(a+1,\frac{r+t}{2}\right)\right|\\ 
& \leq   & 
2 C_{n,\ell,\mu} \Bigg|\left(1+O\left(\frac{1}{t}\right)\right) (t-r)^a \exp \left(-\frac{t}{2}+\frac{r}{2}+O\left(\frac{1}{t^2}\right)\right)\\
&  &
-\left(1+O\left(\frac{1}{t}\right)\right) (r+t)^a \exp \left(-\frac{t}{2}-\frac{r}{2}+O\left(\frac{1}{t^2}\right)\right)\Bigg|  \,. 
\end{eqnarray*}
Hence,  
\begin{eqnarray*} 
\left|F(r,t)- \frac{1}{2 r }\left\{   (r-t) F_0(r-t)+  (r+t) F_0(r+t)\right\} \right|
& \leq   & 
C(r)te^{-t/2} \sum_{k\,odd,\,\,k=-\ell+1}^{\ell-1} t^{\ell- 1-k}t^{a} \\
& \leq   & 
C(r)e^{-t/2}t^{n +\mu - \frac{3}{2} }     \,. 
\end{eqnarray*}
For   odd $\ell$, the arguments are similar. The lemma is proved. \qed
\medskip

The solutions with nontrivial $F_1 $ and their estimate can be easily derived from these formulas by integration over time.\medskip

Although the equation (\ref{EqF}) received the treatment in the textbooks, its application to wave propagation in cosmological backgrounds can still be extended. One must refer to     \cite{Zerilli}, where the author considered the problem of gravitational radiation emitted by a small particle falling into the Schwarzschild field and examined its spectrum in the high-frequency limit.   In \cite{Stewart}, the author  reviewed   equation (\ref{EqF}) and applied it to the computation of the leading terms in the far-field backscatter from outgoing radiation in   Schwarzschild spacetime and of the matter-induced singularities in plane symmetric spacetimes.

\subsection{The Riemann function approach}
\label{SS2.2}

The recursive procedures (\ref{Fellplus1}),(\ref{Felltilde}), and (\ref{FtileF}), can be  reduced to the single formula if we appeal to the Riemann function and modify conditions on the functions at infinity.
Indeed, the equation (\ref{EqF}) can be written in the form 
that leads to the problem for the   unknown function $R (r,t) =rF(r,t) $. In this subsection, we  provide an explicit formula for the solution to the Cauchy problem
\begin{eqnarray}
\label{Z}
\cases{\dsp R_{tt}(r,t)
-   R_{rr} (r,t)+\frac{\ell(\ell+1)}{r^2}R (r,t) =0, \quad r > 0, \quad t \geq 0\cr
{}\cr
R(r,0)=R_0(r),\quad R_t(r,0)=R_1(r), \quad r > 0,\cr
{}\cr
R(0,t)= 0 ,\quad    \quad t\geq 0.}
\end{eqnarray} 
For the Klein-Gordon equation, 
we  set the  initial functions as follows:
\begin{eqnarray}
\label{aroundzeroR}
\hspace{-0.9cm} &  &
R_k(r)=O(r^{\mu +\frac{1}{2}}) \,\, \mbox{\rm as } \,\, r \searrow 0,\quad k=0,1 \,.
\end{eqnarray} 
The functions  $R_0(r)$ and $R_1(r) $  will be continued to negative $r$ as   odd functions for even $\ell$ and as   even functions for odd $\ell$. 

\begin{lemma}
\label{L2.1}
Let $\ell$ be positive integer number and   assume that   $ R_0= R_0(r)$ and $R_1= R_1(r) $  are  odd functions for even $\ell$ and they are  even functions for odd  $\ell$. Assume also $ R_0, R_1 \in C^\infty(\R_+)$, condition (\ref{aroundzeroR}) is fulfilled  and that $\mu  >\ell-\frac{3}{2}$. Then the solution to the problem (\ref{Z}) with the property (\ref{integrable}) is given by 
\begin{eqnarray*} 
&  &
 R (r ,t ) \\
& = &
\!\! \frac{1}{2} [ R_0(r-t )   +  R_0(r+t) ]
-\frac{1}{4 } \ell  (\ell +1) \frac{t}{r}\int_{r-t}^{r+t}  R_0(s)  
\frac{ 1 }{s} F \left(1-\ell ,\ell +2;2;\frac{t^2-(r-s )^2}{4 r s }\right)ds\\
&  &
+ \int_0^t \frac{1}{2} [ R_1(r-\tau )   +  R_1(r+\tau ) ]\,d \tau\\
&  &
-\int_0^t \frac{1}{4 } \ell  (\ell +1) \frac{\tau}{r}\int_{r-\tau}^{r+\tau}  R_1(s)  
\frac{ 1 }{s} F \left(1-\ell ,\ell +2;2;\frac{\tau^2-(r-s )^2}{4 r s }\right)\,ds\,d \tau.
\end{eqnarray*}
 \end{lemma}
 \medskip
 
 \noindent
 {\bf Proof.} 
We skip the proof since it is standard and can be easily reconstructed. (See, e.g., \cite{Stewart,Yag_Galst_CMP} and references therein.)\qed 
\medskip

In particular, for $\ell=5$, in addition  to the cases obtained in subsection~\ref{SS2.1},  we obtain 
 \begin{eqnarray*} 
&  &
F (r ,t ) \\ 
& = &
 \frac{1}{2r} [  (r-t )F_0(r-t )   +  (r+t )F_0(r+t) ]\\
&  &
+\frac{1}{2 r^6}t \left\{\frac{1}{2^7} \left(-315 r^8+1260 r^6 t^2-1890 r^4 t^4+1260 r^2 t^6-315 t^8\right) \int_{r-t}^{r+t}  F _0(\tau)  
 \frac{1}{\tau ^4} \,d\tau \right.\\
&  &
+\frac{1}{2^5} \left(-105 r^6+525 r^4 t^2-735 r^2 t^4+315 t^6\right)\int_{r-t}^{r+t}  F _0(\tau) \frac{1 }{\tau ^2} 
  \,d\tau\\
&  &
+\frac{1}{2^6} \left(-225 r^4+1050 r^2 t^2-945 t^4\right)\int_{r-t}^{r+t}  F _0(\tau)  
 \,d\tau\\
&  &
\left. + \frac{1}{2^5} \left(315 t^2-105 r^2\right) \int_{r-t}^{r+t}  F _0(\tau)  
 \tau ^2\,d\tau 
-\frac{315 }{2^7}\int_{r-t}^{r+t}  F _0(\tau)  
 \tau ^4 \,d\tau\right\}\,.
\end{eqnarray*}
\begin{remark}
The formula of Lemma~\ref{L2.1} as well as the recursive procedures (\ref{Fellplus1}),  (\ref{Felltilde}), and (\ref{FtileF}) exhibit the finite propagation speed property. In particular, for a fixed time, in the fixed bounded domain, the values at infinity of functions $R=R_0(r) $ and $R=R_1(r)$ do  not affect the solution.
\end{remark}

\subsection{Solution by Hankel transform}
\label{SS2.3}

In this subsection we derive one more representation of the solution to the problem for all natural numbers $\ell$ under the assumption of the decay of solution at infinity that is weaker than (\ref{aroundzeroEric_Exp}). 
Consider the problem for the wave equation (\ref{WE}).  
We appeal to the   $\nu$th-order Hankel transform 
\[
F_\nu(s):= \int_0^\infty rf(r)J_\nu(sr)\,dr\,.
\]
with $\nu>-1/2$  (\cite[Sec.11.5.1]{Polyanin-M}, \cite[Sec.14.3]{Watson}, \cite[Sec.14.3]{Zemanian}). Let $f(r) $ be a function defined for $r> 0$.  
The inversion formula yields
\[
f(r) = \int_0^\infty sF_\nu(s)J_\nu(sr)\,ds\,.
\] 
Let $ Y_{\ell m}(\theta, \phi)$, $| m|\leq \ell$, $\ell=0,1,\ldots$, be spherical harmonics. We    skip the proof of the following statement.

\begin{lemma}
\label{L2.3}
Assume that $ F_0(r),  F_1(r) \in    C^1  (0,\infty)$, and 
\[
\int_0^\infty  r ^{3/2}|F_i(r)|\,dr<\infty,\quad    F_i(r)=O(r^{-k}), \quad r \to \infty, \quad k> 2,\quad i=0,1.
\]  
Then the solution  to the problem
\begin{eqnarray*} 
V_{tt}-  \bigtriangleup V =0, \quad V(x,0)= 0, \quad V_t(x,0)= F_1(r)  Y_{\ell m}(\theta, \phi) \,,\quad  t \in {\mathbb R}_+,\, x \in {\mathbb R}^3\,,
\end{eqnarray*}
is 
\begin{eqnarray*}
V (r,\theta, \phi,t)  
& =  &
Y_{\ell m}(\theta, \phi)\frac{1}{\sqrt{r}}\int_0^\infty 
  \left( \int_0^\infty F_1(\rho )  J_{  \ell+\frac{1}{2}}(\lambda  \rho  ) \,\rho ^{3/2} d\rho  \right) \sin (\lambda t)   J_{ \ell+\frac{1}{2}}(r \lambda )\,d\lambda\,  ,   
\end{eqnarray*}
while the solution to the problem
\begin{eqnarray*} 
v_{tt}-  \bigtriangleup v =0, \quad v(x,0)=F_0(r)  Y_{\ell m}(\theta, \phi) 
, \quad v_t(x,0)=0\,,
\end{eqnarray*}
is given by
\begin{eqnarray*}
v (r,\theta, \phi,t)  
& =  &
Y_{\ell m}(\theta, \phi)\frac{1}{\sqrt{r}}\int_{0}^\infty \left(\int_{0}^\infty      F_0(\rho) J_{\ell+\frac{1}{2}}(\rho \lambda )\rho^{3/2}\, d\rho\right)\cos \left( \lambda t \right) J_{\ell+\frac{1}{2}}(r \lambda )\lambda\, d\lambda\,.
\end{eqnarray*}
\end{lemma}
\medskip

In particular, we can set
\begin{eqnarray*}
R ( -r,t) 
& :=  & 
(-1)^{\ell+1}R ( r,t)  \quad \mbox{\rm for all} \quad  r > 0\,.  
\end{eqnarray*}
Thus, the kernel $K(r,t; \rho ) $ of the integral transform 
\begin{eqnarray*}
R (r,t)  
& =  &
\int_0^\infty 
 R_1(\rho ) K(r,t; \rho ) \,\, d\rho  , \quad r > 0, \quad t \geq 0\,,
\end{eqnarray*}
is
\begin{eqnarray*}
K(r,t; \rho ) 
& =  &
  \int_0^\infty   \sin (\lambda t)  \sqrt{\rho }J_{  \ell+\frac{1}{2}}(\lambda \rho ) \sqrt{r}  J_{ \ell+\frac{1}{2}}(r \lambda )\,\,   d\lambda\,,\quad r > 0\,,\quad \rho > 0\,,\quad t \geq 0\,.
\end{eqnarray*}
Note that the function
\begin{eqnarray*}
k(r,t; \rho ,\lambda) 
& =  &
  \sin (\lambda t) \sqrt{\rho } J_{  \ell+\frac{1}{2}}(\lambda \rho ) \sqrt{r}  J_{ \ell+\frac{1}{2}}(r \lambda )\,,\quad r > 0, \,\, t>0,\,\,\lambda>0,\,\,\rho>0\,, 
\end{eqnarray*}
solves the equation 
\begin{eqnarray*}  
&  &
 \dsp k_{tt}(r,t; \rho,\lambda) 
-   k_{rr}(r,t; \rho,\lambda)+\frac{\ell(\ell+1)}{r^2}k (r,t; \rho ,\lambda) =0, \quad t > 0,\quad r > 0\,. 
\end{eqnarray*}

The isometry property in the weighted Sobolev spaces of the Hankel transform is studied in \cite{Airapetyan,Witt}, and it can be extended to the problem of Theorem~\ref{T3.1}.  

\section{Proof of Theorems~\ref{T1.1},\ref{T3.1} }
\label{S3}

\setcounter{equation}{0}

To illustrate ITA, we first present the case of the Minkowski space. Consider in the Minkowski space the spherical Klein-Gordon fields, which, for instance, can be emitted  by the wave function of the pionic atom. 
Application of the ITA    to the Klein-Gordon equation in the Minkowski space leads to 
the following formula  
\begin{eqnarray*} 
u (x,t )
   & =  & 
 v_ {\varphi_0}(x,t)\\
&  &
-  \int_{0 }^{t   }  \frac{m_0 t}{\sqrt{t^2-\tau^2} } J_1\left(m_0\sqrt{t^2-\tau^2} \right)    v_ {\varphi_0} (x,\tau ) \, d\tau  
 +  \int_{0 }^{t   }   J_0\left(m_0\sqrt{t^2-\tau^2} \right)     v_ {\varphi_1} (x,\tau ) \,  d\tau \nonumber     
\end{eqnarray*}
for the solution of the problem
\begin{eqnarray*}
u_{tt}-\Delta u +m_0^2 u =0  \quad \mbox{\rm in}   \quad  {\mathbb
R}^{n+1}, \quad u(x,0 )=\varphi_0 (x), \quad   u_t(x,0 )=\varphi_1 (x) \quad  \mbox{\rm
in}  \quad  {\mathbb R}^{n}\,.
\end{eqnarray*}
Here   $v_ \varphi(x,r ) $ is a solution to the  problem
\begin{eqnarray*}
v_{tt}-\Delta v  =0  \quad \mbox{\rm in}   \quad  {\mathbb
R}^{n+1}, \quad v(x,0 )=\varphi  (x), \quad   v_t(x,0 )=0   \quad  \mbox{\rm
in}  \quad  {\mathbb R}^{n}\,.
\end{eqnarray*}  
One can use any form of the spherical waves given in Section~\ref{SWM}. We choose  Lemma~\ref{L2.3} and  
 skip an elementary  proof of the following statement.

\begin{theorem}
Assume that $ F_0(r),  F_1(r) \in    C^1  (0,\infty)$, and 
\[
\int_0^\infty  r ^{3/2}|F_i(r)|\,dr<\infty,\quad    F_i(r)=O(r^{-k}), \quad r \to \infty, \quad k> 2,\quad i=0,1.
\] 
The  spherical   Klein-Gordon field in the Minkowski space is given by 
\begin{eqnarray*}
\Phi (r,\theta, \phi,t) 
& = & 
 Y_{\ell m}(\theta, \phi) \frac{1}{r}\Bigg\{\int_0^\infty 
  \left( \int_0^\infty F_0(\rho)  J_{  \ell+\frac{1}{2}}(\lambda \rho ) \,\rho^2 d\rho \right) \cos (\lambda t)    J_{ \ell+\frac{1}{2}}(r \lambda )\,\lambda d\lambda  \\
&  &
-   \int_{0 }^{t   }  \frac{mt}{\sqrt{t^2-\tau^2} } J_1\left(m\sqrt{t^2-\tau^2} \right)   \int_0^\infty 
  \left( \int_0^\infty F_0(\rho)  J_{  \ell+\frac{1}{2}}(\lambda \rho ) \,\rho^2 d\rho \right) \cos (\lambda \tau)    J_{ \ell+\frac{1}{2}}(r \lambda )\,\lambda d\lambda  \, d\tau    \\
 &    & 
 +  \int_{0 }^{t   }   J_0\left(m\sqrt{t^2-\tau^2} \right)\int_0^\infty 
  \left( \int_0^\infty F_1(\rho)  J_{  \ell+\frac{1}{2}}(\lambda \rho ) \,\rho^2 d\rho \right) \cos (\lambda \tau)    J_{ \ell+\frac{1}{2}}(r \lambda )\,\lambda d\lambda  \,  d\tau   \Bigg\}
\end{eqnarray*} 
 for all $ x \in {\mathbb R}^3$,   $ t>0 $.
\end{theorem}
\medskip

Next we turn to the de~Sitter universe. 
The ITA applied to the initial value  problem (\ref{1.10n}) (written in the Cartesian coordinates $x \in {\mathbb R}^n$),  for the equation  
\begin{equation}
\label{psieq} 
    \frac{\partial^2 \psi }{\partial t^2}
+ nH        \frac{\partial \psi }{\partial t}
-e^{-2Ht}{\mathcal{A}}(x,\partial_x)\psi +   \frac{m_n^2}{\hbar^2} \psi=0\,,
\end{equation} 
with the operator $ {\mathcal{A}}(x,\partial_x)= \sum _{|\alpha|\leq 2 } a_\alpha(x)   \partial_x^\alpha$, leads to the  formula for the solution
\begin{eqnarray}
\label{psi} 
&   &
\psi (x,t) \\
&  =  &
 e ^{-\frac{n-1}{2}Ht} v_{\psi_0}  (x, \phi (t))  
+ \, e ^{-\frac{n}{2}Ht}\int_{ 0}^{\phi (t)}  \left[ 2K_0( s,t;M)    
+   nHK_1( s,t;M)  \right] v_{\psi_0 } (x,  s)
  ds  \nonumber \\
&  &
+\, 2e ^{-\frac{n}{2}Ht}\int_{0}^{\phi (t) }  v_{\psi_1 } (x,  s)
  K_1( s,t;M)   ds,\quad  \forall x \in \Omega \subseteq \R^n \,,\,\, \forall t \in I=[0,T] \subseteq [0,\infty), \nonumber 
\end{eqnarray}
of the Cauchy problem for the equation (\ref{psieq}) with the initial data $\psi (x,0) = \psi_0 (x)$ and $\psi_t (x,0) = \psi_1 (x)$. 
Here the kernels  $K_0( s,t;M) $  and $K_1( s,t;M) $,   are defined in  (\ref{K0M}) and (\ref{K1M}), (see \cite{Yag_Galst_CMP,JMP_2024}), 
 respectively, 
where $0<T\leq\infty$, $\phi (t):= (1-e^{-Ht} )/H$, $M^2 =  \frac{n^2 H^2}{4 } -\frac{m^2c^4}{\hbar^2} \in {\mathbb C}$, and   the function \, $v=  v_\varphi(x, t)  $ \, is a   solution   of the   problem
\begin{equation}
\label{vphi}
\cases{
 v_{tt}-   {\mathcal A}(x,\partial_x)  v =0, \quad \forall x \in \Omega \,,\quad \forall t \in [0,(1-e^{-HT})/H]\,, \cr
 v(x,0)= \varphi (x), \quad v_t(x,0)=0\,,\quad \forall x \in \Omega\,.  
} 
\end{equation} 
   Here, 
the function $v_\varphi  (x, \phi (t) )$  coincides with the value $v(x, \phi (t) ) $
of the solution $v(x,t)$ of the Cauchy problem (\ref{vphi}) with $ {\mathcal A}(x,\partial_x)=\Delta$. In the case of $n=3$ and $m_n^2   =2H^2\hbar^2 $ the field is huygensian, that is, it obeys the Huygens' principle (see \cite{JMP2013}) and the kernels are simplified as follows:
\[     
K_0\left(r,t;\frac{1}{2}H\right)
   =   -\frac{1}{4} H e^{\frac{H t}{2}}  \,,  \quad 
K_1\left(r,t;\frac{1}{2}H\right)
  =   \frac{1}{2} e^{\frac{H t}{2}} \,.  
\]
 (For the definition of the Huygens' principle see, e.g.,  \cite{Gunther}.)
\medskip

Thus, one can regard ITA as an analytical mechanism that, from the massless virtual field in the Minkowski space, generates a massive field in the curved spacetime.   The integral transform condenses the evolution of that virtual scalar wave in the conformal time from its entry into the expanding universe to the moment of observation.

\subsection{Proof of Theorem~\ref{T1.1}}

 From  (\ref{psi}) we derive
\begin{eqnarray} 
\label{3.19} 
\Phi (r,\theta, \phi,t) \nonumber   
& = &
e^{-Ht} v_{F_0 Y_{\ell m}}  \left(x, \phi (t) \right) 
+ \,  e^{-H\frac{3}{2}t}\int_{ 0}^{\phi (t)} v_{F_0 Y_{\ell m}}  \left(x, \tau\right) \left(2  K_0 \left(\tau ,t \right)+ 3K_1\left(\tau ,t \right)\right)\,  d\tau   \nonumber \\
& &
+\, 2e^{-\frac{3}{2}Ht}\int_{0}^{\phi (t)}   v_{F_1 Y_{\ell m} } \left(x, \tau \right)
  K_1\left(\tau,t \right)\, d\tau\,.
\end{eqnarray}
Then we use the formula 
from Lemma~\ref{L2.1} with $F_1=0$, that is,
\begin{eqnarray}
\label{3.21}
v_{F_k Y_{\ell m}}  \left(x, t \right) 
& = &
 Y_{\ell m}(\theta, \phi)\Bigg\{ \frac{1}{2r} [ (r-t )F_k(r-t )   +  (r+t)F_k(r+t) ]\\
&  &
-\frac{1}{4 } \ell  (\ell +1) \frac{t}{r^2}\int_{r-t}^{r+t}  F_k(s)  
F \left(1-\ell ,\ell +2;2;\frac{t^2-(r-s )^2}{4 r s }\right)ds\Bigg\},\quad k=0,1\,. \nonumber
\end{eqnarray}
Hence, with $\phi (t):= \frac{1-e^{-Ht}}{H} $ and $\tau:=\phi (t)s =\frac{1-e^{-Ht}}{H}s$ we obtain
\begin{eqnarray*} 
&  &
\Phi (r,\theta, \phi,t) \\ 
& = &
e^{-Ht} Y_{\ell m}(\theta, \phi)\Bigg\{ \frac{1}{2r} [ (r-\phi (t) )F_0(r-\phi (t) )   +  (r+\phi (t))F_0(r+\phi (t)) ]\\
&  &
-\frac{1}{4 } \ell  (\ell +1) \frac{\phi (t)}{r^2}\int_{r-\phi (t)}^{r+\phi (t)}  F_0(s)  
F \left(1-\ell ,\ell +2;2;\frac{\phi (t)^2-(r-s )^2}{4 r s }\right)ds\Bigg\} \\
&  &
+ \,  e^{-H\frac{3}{2}t}Y_{\ell m}(\theta, \phi)\int_{ 0}^{\phi (t)}  \Bigg\{ \frac{1}{2r} [ (r-\tau )F_0(r-\tau )   +  (r+\tau)F_0(r+\tau) ]\\
&  &
-\frac{1}{4 } \ell  (\ell +1) \frac{\tau}{r^2}\int_{r-\tau}^{r+\tau}  F_0(s)  
F \left(1-\ell ,\ell +2;2;\frac{t^2-(r-s )^2}{4 r s }\right)ds\Bigg\} \left(2  K_0 \left(\tau ,t \right)+ 3K_1\left(\tau ,t \right)\right)\,  d\tau  \nonumber \\
& &
+\, 2e^{-\frac{3}{2}Ht}Y_{\ell m}(\theta, \phi)\int_{0}^{\phi (t)}   \Bigg\{ \frac{1}{2r} [ (r-\tau )F_1(r-\tau )   +  (r+\tau)F_1(r+\tau) ]\\
&  &
-\frac{1}{4 } \ell  (\ell +1) \frac{\tau}{r^2}\int_{r-\tau}^{r+\tau}  F_1(s)  
F \left(1-\ell ,\ell +2;2;\frac{\tau^2-(r-s )^2}{4 r s }\right)ds\Bigg\}
  K_1\left(\tau,t \right)\, d\tau\,.
\end{eqnarray*}

For the huygensian field  with   
$
\varphi _0(x)=\Phi (x,0)$ and $\varphi _1(x)=\Phi _t(x,0) 
$
of   (\ref{1.10n}),
   the field in the de~Sitter FLRW spacetime is (see \cite{JMP2013})  
\begin{eqnarray}
\label{Huygen} 
\Phi  (r,\theta, \phi,t) 
& = &
e^{-Ht} v_{\varphi_0}  \left( x, \phi (t) \right)
+ \,   e^{-Ht}    V_{\varphi_0+\varphi_1}  \left( x, \phi (t) \right)  \,,
\end{eqnarray}
where $ V_{\varphi_0+\varphi_1}   ( x, t  )$ is a solution of the problem
\begin{eqnarray*}  
V_{tt}-  \Delta V =0, \quad V(x,0)=0 , \quad V_t(x,0)=\varphi_0(x)+\varphi_1(x)\,,
\end{eqnarray*}
and $ V_{\varphi_0+\varphi_1}(x,t)=\int_0^t v(x,\tau)\,d \tau$. 
Next, we use (\ref{3.21}),
\begin{eqnarray*}
v_{F_0 Y_{\ell m}}  \left(x, \phi (t) \right) 
& = &
 Y_{\ell m}(\theta, \phi)\Bigg\{ \frac{1}{2r} [ (r-\phi (t) )F_0(r-\phi (t) )   +  (r+\phi (t))F_0(r+\phi (t)) ]\\
&  &
-\frac{1}{4 } \ell  (\ell +1) \frac{\phi (t)}{r^2}\int_{r-\phi (t)}^{r+\phi (t)}  F_0(s)  
F \left(1-\ell ,\ell +2;2;\frac{\phi (t)^2-(r-s )^2}{4 r s }\right)ds\Bigg\},  
\end{eqnarray*} 
and
\begin{eqnarray*}
&  &
V_{(F_0+F_1) Y_{\ell m}}  \left(x, \phi (t) \right) \\
& = &
 Y_{\ell m}(\theta, \phi)\int_0^{\phi (t)} \Bigg\{ \frac{1}{2r} \Big( (r-\tau )[F_0(r-\tau ) +F_1(r-\tau )]  +  (r+\tau)[F_0(r+\tau)+F_1(r+\tau)] \Big)\\
&  &
-\frac{1}{4 } \ell  (\ell +1) \frac{\tau}{r^2}\int_{r-\tau}^{r+\tau}  (F_0(s)+F_1(s))   
F \left(1-\ell ,\ell +2;2;\frac{\tau^2-(r-s )^2}{4 r s }\right)ds\Bigg\}d \tau,  
\end{eqnarray*}
 to obtain the last statement of the theorem.  
The theorem is proved. \qed

\subsection{Proof of Theorem~\ref{T3.1}} 

\noindent
{\bf Proof of Theorem~\ref{T3.1}.} 
If we use the formula  (\ref{3.19}), then   we derive
\begin{eqnarray*} 
&  &
\Phi (r,\theta, \phi,t) \\  
& = &
e^{-Ht} Y_{\ell m}(\theta, \phi)\frac{1}{\sqrt{r}}\int_{0}^\infty \left(\int_{0}^\infty      F_0(\rho) J_{\ell+\frac{1}{2}}(\rho \lambda )\rho^{3/2}\, d\rho\right)\cos \left( \lambda \phi (t)  \right) J_{\ell+\frac{1}{2}}(r \lambda )\lambda\, d\lambda \\
& &
+ \,  e^{-H\frac{3}{2}t}Y_{\ell m}(\theta, \phi)\frac{1}{\sqrt{r}}\int_{ 0}^{1}\int_{0}^\infty \left(\int_{0}^\infty      F_0(\rho) J_{\ell+\frac{1}{2}}(\rho \lambda )\rho^{3/2}\, d\rho\right)\cos \left( \lambda \phi (t)  s \right)  \\
&  &
\times J_{\ell+\frac{1}{2}}(r \lambda )\lambda\, d\lambda\left(2  K_0 \left(\phi (t) s,t \right)+ 3K_1\left(\phi (t) s,t \right)\right)\phi (t) \,  ds  \nonumber \\
& &
+\, 2e^{-\frac{3}{2}Ht}Y_{\ell m}(\theta, \phi)\frac{1}{\sqrt{r}}\int_{0}^1 \int_{0}^\infty \left(\int_{0}^\infty      F_1(\rho) J_{\ell+\frac{1}{2}}(\rho \lambda )\rho^{3/2}\, d\rho\right)\cos \left( \lambda \phi (t)  s \right)  \\
&  &
\times 
J_{\ell+\frac{1}{2}}(r \lambda )\lambda\, d\lambda  \, K_1\left(\phi (t) s,t \right)\phi (t) \, ds \quad \mbox{\rm for all}  \quad t>0\, . 
\end{eqnarray*}
For the huygensian field   we use (\ref{Huygen}) and  the formula from Lemma~\ref{L2.3}, to obtain 
\begin{eqnarray*}
&  &
\Phi  (r,\theta, \phi,t)\\
  & = &
e^{-Ht}  Y_{\ell m}(\theta, \phi)\frac{1}{\sqrt{r}}\Bigg\{ \int_{0}^\infty \left(\int_{0}^\infty      F_0(\rho) J_{\ell+\frac{1}{2}}(\rho \lambda )\rho^{3/2}\, d\rho\right)\cos \left( \lambda \phi (t)  \right) J_{\ell+\frac{1}{2}}(r \lambda )\lambda\, d\lambda \\
&  &
+ 
\int_0^\infty 
  \left( \int_0^\infty (   F_0(\rho)+F_1(\rho))  J_{  \ell+\frac{1}{2}}(\lambda  \rho  ) \,\rho ^{3/2} d\rho  \right) \sin \left(\lambda \phi (t)  \right)   J_{ \ell+\frac{1}{2}}(r \lambda )\,d\lambda \Bigg\}\,. 
\end{eqnarray*}
Theorem~\ref{T3.1} is proved.
\qed
\medskip

\begin{corollary}
\label{C3.2}
Let $ F_1(r)\equiv 0$ and $ F_0(r)$ be of the form (\ref{aroundzeroEric_Exp}). 
If 
$\dsp  3H\hbar  /2\leq  m_n  $, then
\begin{eqnarray*} 
\left| \Phi  (r,\theta, \phi,t)- e ^{- Ht}  v_{F_0 Y_{\ell m}}  (r,\theta, \phi, \phi (t)) \right| 
  & \leq  &
   C(r)e ^{- Ht}\left| Y_{\ell m}(\theta, \phi)\right|(1+t)^{1-\sgn (|M|)}\,,
  \end{eqnarray*}
where $ \phi (t)=(1-e^{-Ht} )/H$. If 
$\dsp  3H\hbar  /2 > m_n  $, then
\begin{eqnarray*}
&  &
\left| \Phi  (r,\theta, \phi,t)- e ^{- Ht}  v_{F_0 Y_{\ell m}}  (r,\theta, \phi, \phi (t)) \right| \\
  & \leq  &
 C(r)e ^{-\frac{3}{2}Ht}\left| Y_{\ell m}(\theta, \phi) \right|  \cases{   e^ {  \frac{1}{2}Ht}  \quad \mbox{\rm if} \quad  3H\hbar  /2>  m_n   \geq \sqrt{2}H \hbar\,, \cr
 e^ { t\sqrt{\frac{9}{4}H^2 -\frac{m_n^2 }{\hbar^2} } }    \quad \mbox{\rm if} \quad m_n  <\sqrt{2}H \hbar\,.}  
\end{eqnarray*}
Here $ v_{F_0 Y_{\ell m}}  (x, \phi (t))$ is given by Lemma~\ref{L2.3},
\begin{eqnarray*}
&  &
v_{F_0 Y_{\ell m}}  (r,\theta, \phi, \phi (t))\\  
&  =&
Y_{\ell m}(\theta, \phi)\frac{1}{\sqrt{r}}\int_{0}^\infty \left(\int_{0}^\infty      F_0(\rho) J_{\ell+\frac{1}{2}}(\rho \lambda )\rho^{3/2}\, d\rho\right)\cos \left( \lambda \phi (t) \right) J_{\ell+\frac{1}{2}}(r \lambda )\lambda\, d\lambda\,.
\end{eqnarray*}
\end{corollary}
\medskip

\noindent
{\bf Proof.} We appeal to the representation (\ref{psi}) with the dimension of   spatial variable $n=3$,
\[
\psi (x,t) 
=
 e ^{- Ht} v_{\psi_0}  (x, \phi (t))   
+ \, e ^{-\frac{3}{2}Ht}\int_{ 0}^{\phi (t)}  \left[ 2K_0( s,t;M)    
+   3HK_1( s,t;M)  \right] v_{\psi_0 } (x,  s)  ds\,,
\]
and to Lemma~\ref{LTail}. 
For the tail, this implies
\begin{eqnarray*}
&  &
\left| \Phi  (r,\theta, \phi,t)- e ^{- Ht}  v_{F_0 Y_{\ell m}}  (x, \phi (t)) \right| \\
  & = &
\left| e ^{-\frac{3}{2}Ht}\int_{ 0}^{\phi (t)}  \left[ 2K_0( s,t;M)    
+   3HK_1( s,t;M)  \right]  v_{F_0 Y_{\ell m}} (x,  s)
  ds  \right| \\
  & \leq &
e ^{-\frac{3}{2}Ht}\left| Y_{\ell m}(\theta, \phi) \right|\int_{ 0}^{\phi (t)}  \left| 2K_0( s,t;M)    
+   3HK_1( s,t;M)  \right| \\
&  &
\times\left( \frac{1}{2 r }\left|     (r-s) F_0(r-s)+  (r+s) F_0(r+s)\right|  +    C(r)e^{-s/2}s^{n +\mu -\frac{3}{2} }  \right) ds\\
  & \leq &
 C(r)e ^{-\frac{3}{2}Ht}\left| Y_{\ell m}(\theta, \phi) \right|\int_{ 0}^{\phi (t)}  \left| 2K_0( s,t;M)    
+   3HK_1( s,t;M)  \right| \\
&  &
\times\Bigg( \frac{1}{2 r }\Bigg[        |r-s|^{\mu +\frac{1}{2} } e^{-|r-s|/2 }|L_{n-\ell-1}^{2\mu }(|r-s|)|\\
&  &
+   |r+s|^{\mu +\frac{1}{2} } e^{-|r+s|/2 }|L_{n-\ell-1}^{2\mu }(|r+s|)|\Bigg] 
+   e^{-s/2}s^{n +\mu -\frac{3}{2} }  \Bigg) ds\\
  & \leq &
 C(r)e ^{-\frac{3}{2}Ht}\left| Y_{\ell m}(\theta, \phi) \right|\int_{ 0}^{\phi (t)}  \left| 2K_0( s,t;M)    
+   3HK_1( s,t;M)  \right| \\
&  &
\times\Bigg(       (r+s)^{\mu +\frac{1}{2}+ n-\ell-1}  
+   e^{-s/2}s^{n +\mu -\frac{3}{2} }  \Bigg) ds\\
  & \leq &
 C(r)e ^{-\frac{3}{2}Ht}\left| Y_{\ell m}(\theta, \phi) \right|\int_{ 0}^{\phi (t)}  \left| 2K_0( s,t;M)    
+   3HK_1( s,t;M)  \right|   ds,\quad \ell=1,2,\ldots,
\end{eqnarray*}
where $v_{F_0 Y_{\ell m}}  (x, \phi (t)) $ can be   written in three different explicit forms of Subsections~\ref{SS2.1}-\ref{SS2.3}.
Recall, $\dsp M :=\left(\frac{9}{4}H^2 -\frac{m_n^2 }{\hbar^2} \right)^{1/2}$; hence, either $M\geq 0$ or $\Re M=0$. 
According to Lemmas~4.2,~4.3~\cite{JMP_2024}, 
if $M\geq 0$ and $H $ is a positive number, then
\begin{eqnarray*}
\int_{ 0}^{\phi (t)} | K_0( s,t;-iM) |   ds 
 & \leq & 
C_{M,H}(1+t)^{1-\sgn (M)} e^{- Ht/2} (e^{Ht}-1) \quad \mbox{for all}\quad t \in [0,\infty)\,,\\
 \int_{0}^{\phi (t) }  | K_1( r,t;-iM) |  dr 
 & \leq &  
   C_{M } (1+t)^{1-\sgn (M)}\frac{1}{ H} (e^{H t}  - 1) (e^{H t}  + 1)^{-1 } \quad \mbox{for all} \quad t \in [0,\infty) \,.
\end{eqnarray*}
According to Lemmas~5.1,~5.2~\cite{JMP_2024}, 
if $  a>-1 $, $  M>0$,  and $\phi (t)= (1-e^{-Ht})/H$, then for large  $t$,
\begin{eqnarray*}
 \int_{ 0}^{1} \phi (t)^{a} s^{a}
\big|  K_1(\phi (t)s,t;M)\big|   \phi (t)\,  ds  
  & \leq & 
 C_{M,a}   e^{  M   t}  \,,\\
 \int_{ 0}^{1} \phi (t)^{a} s^{a}
\big|  K_0(\phi (t)s,t;M)\big|   \phi (t)\,  ds  
  & \leq &   C_{M,a}  
\times \cases{   e^ {  \frac{1}{2}Ht}  \quad \mbox{\rm if} \quad   M<H/2\,, \cr
 e^ { (  M- H) t}    \quad \mbox{\rm if} \quad   M> H/2\,.}  
\end{eqnarray*}
Consequently, if 
$\dsp  3H\hbar  /2> m_n  $, then
\begin{eqnarray*}
&  &
\left| \Phi  (r,\theta, \phi,t)- e ^{- Ht}  v_{F_0 Y_{\ell m}}  (x, \phi (t)) \right| \\
& \leq &
 C(r)e ^{-\frac{3}{2}Ht}\left| Y_{\ell m}(\theta, \phi) \right|\\
&  &
\times \left( e^{\left(\frac{9}{4}H^2 -\frac{m_n^2 }{\hbar^2} \right)^{1/2}   t}  + \cases{   e^ {  \frac{1}{2}Ht}  \quad \mbox{\rm if} \quad  3H\hbar  /2>  m_n   >\sqrt{2}H \hbar\,, \cr
 e^ { (\left(\frac{9}{4}H^2 -\frac{m_n^2 }{\hbar^2} \right)^{1/2}- H) t}    \quad \mbox{\rm if} \quad m_n  <\sqrt{2}H \hbar\,,} \right) \\
& \leq &
 C(r)e ^{-\frac{3}{2}Ht}\left| Y_{\ell m}(\theta, \phi) \right|  \cases{   e^ {  \frac{1}{2}Ht}  \quad \mbox{\rm if} \quad  3H\hbar  /2>  m_n   >\sqrt{2}H \hbar\,, \cr
 e^ { t\sqrt{\frac{9}{4}H^2 -\frac{m_n^2 }{\hbar^2} } }    \quad \mbox{\rm if} \quad m_n  <\sqrt{2}H \hbar\,.}    
\end{eqnarray*}
If 
$\dsp  3H\hbar  /2 \leq m_n  $, then
\begin{eqnarray*}
&  &
\left| \Phi  (r,\theta, \phi,t)- e ^{- Ht}  v_{F_0 Y_{\ell m}}  (x, \phi (t)) \right| \\
  & \leq  &
 C(r)e ^{-\frac{3}{2}Ht}\left| Y_{\ell m}(\theta, \phi) \right|\int_{ 0}^{\phi (t)}  \left| 2K_0( s,t;M)    
+   3HK_1( s,t;M)  \right|   ds \\
& \leq &
 C(r)e ^{- Ht}\left| Y_{\ell m}(\theta, \phi) \right|(1+t)^{1-\sgn (|M|)}\,.
\end{eqnarray*}
For the case of huygensian field $m_n  =\sqrt{2}H \hbar$, the solution (\ref{Huygen}) implies
\begin{eqnarray*}
&  &
\left| \Phi  (r,\theta, \phi,t)- e ^{- Ht} v_{\psi_0}  (x, \phi (t)) \right| \\
&  \leq &
 e ^{- Ht}H\int_{ 0}^{\phi (t)} \left|   v_{\psi_0 } (x,  s)  ds  \right|  \\ 
& \leq &
e ^{- Ht}H\left| Y_{\ell m}(\theta, \phi) \right|\\
&  &
\times \int_{ 0}^{\phi (t)} \left|     \frac{1}{2 r }\left|     (r-s) F_0(r-s)+  (r+s) F_0(r+s)\right|  +    C(r)e^{-s/2}s^{n +\mu -\frac{3}{2}}    \right| ds   \\ 
& \leq &
C(r)e ^{- Ht} \left| Y_{\ell m}(\theta, \phi) \right|\,.
\end{eqnarray*}
The corollary is proved.
\qed
\medskip

For the  pionic atom, the modification is evident since $F_1(r)=-\frac{i}{\hbar}E_n F_0(r) $. If $\ell=0$, then the explicit formula of Subsection~\ref{SS2.1} can be used to estimate a solution. 
In particular, for the ground state $n=1$ and $\ell=0$, in the case of huygensian field $m_n  =\sqrt{2}H \hbar$, according to (\ref{Huygen}), at $r=1/H$, we obtain
\begin{eqnarray*} 
\Phi  \left(\frac{1}{H},\theta, \phi,t\right) 
& = &
e^{-Ht} H^{\frac{1}{2} -\mu }\frac{1}{4\sqrt{\pi}}   
\Bigg[\Bigg\{    e^{-H(\mu +\frac{1}{2}) t-   \frac{ 1}{2H}e^{-Ht}   }  
+   \left(   2-e^{-Ht} \right)^{\mu +\frac{1}{2} } e^{-   \frac{1}{2H}(2-e^{-Ht})   }   \Bigg\}\\
&  &
+     \int_0^t \Bigg\{     e^{-H(\mu +\frac{1}{2})s-   \frac{1}{2H} e^{-Hs}   }  
+   \left(   2-e^{-Hs} \right)^{\mu +\frac{1}{2} } e^{-   \frac{1}{2H} (2-e^{-Hs})  }   \Bigg\}\,ds  \Bigg] \,.
\end{eqnarray*}
This indicates the rate of decay of the $L^\infty$-norm of solution for a large $t$.

\section*{Acknowledgments}

The author is indebted to the reviewers  for the valuable
remarks and suggestions, which improved the text.

\end{document}